\documentclass[11pt,reqno]{amsart}

\theoremstyle{plain}
\newtheorem{theorem} {Theorem}[section]
\newtheorem{lemma}[theorem] {Lemma}

\newtheorem{corollary}[theorem] {Corollary}

\theoremstyle{definition}
\newtheorem{definition}[theorem] {Definition}

\theoremstyle{remark}
\newtheorem{remark}[theorem] {Remark}

\numberwithin{equation}{section}

\usepackage{amsfonts}
\usepackage{amssymb}
\usepackage{amscd}

\newcommand{\R}{{\mathbb R}}
\newcommand{\Z}{{\mathbb Z}}
\newcommand{\N}{{\mathbb N}}

\newcommand{\PP}{{\mathcal P}}

\newcommand{\CC}{{\mathbb C}}

\newcommand{\TT}{{\mathcal T}}

\newcommand{\E}{{\mathcal E}}

\newcommand{\BB}{{\mathfrak B}}
\newcommand{\al}{{\alpha}}
\newcommand{\la}{{\lambda}}
\newcommand{\sa}{{\sigma}}

\newcommand{\iy}{{\infty}}

\newcommand{\vep}{{\varepsilon}}

\newcommand{\de}{{\delta}}

\newcommand{\be}{{\beta}}

\newcommand{\bna}{\begin{eqnarray}}
\newcommand{\ena}{\end{eqnarray}}
\newcommand{\ba}{\begin{eqnarray*}}
\newcommand{\ea}{\end{eqnarray*}}
\newcommand{\beq}{\begin{equation}}
\newcommand{\eeq}{\end{equation}}

\begin{document}

\title[Constants in Multivariate Inequalities]
{Sharp Constants of Approximation Theory. III. Certain Polynomial Inequalities
of Different Metrics on Convex Sets}
\author{Michael I. Ganzburg}
 \address{Department of Mathematics\\ Hampton University\\ Hampton,
 VA 23668\\USA}
 \email{michael.ganzburg@hamptonu.edu}
 \keywords{Sharp constants, multivariate Markov-Bernstein-Nikolskii
  type inequality, algebraic polynomials,
 entire functions of exponential type.}
 \subjclass[2010]{Primary 41A17, 41A63, Secondary 26D10}

 \begin{abstract}
 Let $V\subset\R^m$ be a centrally symmetric convex body and let $V^*\subset\R^m$ be its polar.
 We prove
limit relations between the sharp constants in the
multivariate
Markov-Bernstein-Nikolskii
type inequalities
 for algebraic polynomials on $V^*$ and the corresponding constants for
 entire functions of exponential type with the spectrum in $V$.
 \end{abstract}
 \maketitle

 \section{Introduction}\label{S1}
\setcounter{equation}{0}
\noindent
We continue the study of the sharp constants in multivariate inequalities
of approximation theory
that began in \cite{G2018} and \cite{G2019}. In this paper we prove
limit relations between the sharp constants in the
multivariate
 Markov-Bernstein-Nikolskii
type inequalities for algebraic polynomials on convex sets
 and entire functions of exponential type.
\vspace{.12in}\\
\textbf{Notation.}
Let $\R^m$ be the Euclidean $m$-dimensional space with elements
$x=(x_1,\ldots,x_m),\, y=(y_1,\ldots,y_m),
\,t=(t_1,\ldots,t_m),\,u=(u_1,\ldots,u_m)$,
the inner product $t\cdot y:=\sum_{j=1}^mt_jy_j$,
and the norm $\vert t\vert:=\sqrt{t\cdot t}$.
Next, $\CC^m:=\R^m+i\R^m$ is the $m$-dimensional complex space with elements
$z=(z_1,\ldots, z_m)=x+iy$
and the norm $\vert z\vert:=\sqrt{\vert x\vert^2+\vert y\vert^2}$;
$\Z^m$ denotes the set of all integral lattice points in $\R^m$;
and $\Z^m_+$ is a subset of $\Z^m$
of all points with nonnegative coordinates.
We also use multi-indices $\be=(\be_1,\ldots,\be_m)\in \Z^m_+$
and $\al=(\al_1,\ldots,\al_m)\in \Z^m_+$
with
 \ba
 \vert\be\vert:=\sum_{j=1}^m\be_j,\qquad
 \vert\al\vert:=\sum_{j=1}^m\al_j,\qquad
 y^\be:=y_1^{\be_1}\cdot\cdot\cdot y_m^{\be_m}, \qquad
 D^\al:=\frac{\partial^{\al_1}}{\partial y_1^{\al_1}}\cdot\cdot\cdot
 \frac{\partial^{\al_m}}{\partial y_m^{\al_m}}.
 \ea
 \noindent
 Let  $\PP_{n,m}$ be a set of all
  polynomials $P(y)=\sum_{\vert\be\vert\le n}c_\be y^\be$
   in $m$ variables of total degree at most $n,\,n\in Z^1_+$,
   with complex coefficients.
Given $\sa\in\R^m,\,\sa_j>0,\,1\le j\le m$, and $M>0$, let
$\Pi^m(\sa):=\{t\in\R^m: \vert t_j\vert\le\sa_j, 1\le j\le m\},\,
Q^m(M):=\{t\in\R^m: \vert t_j\vert\le M, 1\le j\le m\}$,
and
$\BB^m(M):=\{t\in\R^m: \vert t\vert\le M\}$
be the $m$-dimensional parallelepiped, cube, and ball, respectively.

Let $L_r(E)$ be the space of all measurable complex-valued functions $F$
 on a measurable set $E\subseteq\R^m$  with the finite quasinorm
 \ba
 \|F\|_{L_r(E)}:=\left\{\begin{array}{ll}
 \left(\int_E\vert F(x)\vert^r dx\right)^{1/r}, & 0<r<\iy,\\
 \mbox{ess} \sup_{x\in E} \vert F(x)\vert, &r=\iy.
 \end{array}\right.
 \ea
 This quasinorm allows the following "triangle" inequality:
 \beq\label{E1.1}
 \left\|\sum_{j=1}^s F_j\right\|^{\tilde{r}}_{L_r(E)}
 \le \sum_{j=1}^s \left\|F_j\right\|^{\tilde{r}}_{L_r(E)},
 \qquad F_j\in L_r(E),\qquad
 1\le j\le s,
 \eeq
 where $s\in\N:=\{1,\,2,\ldots\}$ and $\tilde{r}:=\min\{1,r\}$ for $r\in(0,\iy]$.

 In this paper we will need certain definitions and properties of
 convex bodies in $\R^m$.
 We first define the \emph{width} $w(K)$ of a convex body $K$ in $\R^m$ as
 the minimal distance between two parallel supporting hyperplanes of $K$.

 Next, throughout the paper $V$ is a centrally symmetric (with respect to the origin)
 closed
 convex body in $\R^m$ and
 $V^*:=\{y\in\R^m: \forall\, t\in V, \vert t\cdot y\vert \le 1\}$
 is the \emph{polar} of $V$.
 It is well known that $V^*$ is a centrally symmetric (with respect to the origin)
 closed
 convex body in $\R^m$ and $V^{**} =V$ (see, e.g., \cite[Sect. 14]{R1970}).
 The set $V$ generates the following dual norms
 on $\R^m$ and $\CC^m$ by
 \ba
 \|y\|_V^*:=\sup_{t\in V}\vert t\cdot y\vert,\quad y\in\R^m;\qquad
 \|z\|_V^*:=\sup_{t\in V}\left\vert\sum_{j=1}^m t_jz_j\right\vert,\quad z\in\CC^m.
 \ea
 Note also that $V^*$ is the unit ball in the norm $\|\cdot\|_V^*$ on $\R^m$
 with the boundary $\partial(V^*)$, the width $w(V^*)$, and
 the $m$-dimensional volume $\vert V^*\vert_m$.
 For example, if $\sa\in\R^m,\,\sa_j>0,\,1\le j\le m$, and $V=\left\{t\in\R^m:
 \left(\sum_{j=1}^m\vert t_j/\sa_j\vert^{\mu}\right)^{1/\mu}\le 1\right\}$, then
 for $y\in\R^m,\,
 \|y\|_V^*=\left(\sum_{j=1}^m\vert \sa_j y_j\vert^\la\right)^{1/\la}$
 and $V^*=\left\{y\in\R^m:
 \left(\sum_{j=1}^m\vert \sa_j y_j\vert^\la\right)^{1/\la}\le 1\right\}$,
 where $\mu\in[1,\iy],\,\la\in[1,\iy]$, and $1/\mu+1/\la=1$.
 In particular,
 $\|y\|_{\Pi^m(\sa)}^*=\sum_{j=1}^m\sa_j\vert y_j\vert,\,
 \|y\|_{Q^m(M)}^*=M\sum_{j=1}^m\vert y_j\vert$,
 and $\|y\|_{\BB^m(M)}^*=M\vert y\vert$.

 Let $l_y$ be a straight line passing through
   a fixed point $y\in\R^m$ and the origin,  and
   let
   $v_y\in l_y\cap\partial(V^*)$
   be the closest point to $y$. We also need
   the equivalent definition of the dual norm in $\R^m$
   given by the formula
   (see, e.g., \cite[Theorem 14.5]{R1970})
   \beq\label{E1.1a}
   \|y\|_V^*=\vert y\vert/\vert v_y\vert.
   \eeq

 Given $a>0$, the set of all trigonometric polynomials
 $T(x)=\sum_{k\in aV\cap \Z^m}c_ke^{ik\cdot x}$ with complex
 coefficients is denoted by $\TT_{aV}$.

  \begin{definition}\label{D1.1}
 We say that an entire function $f:\CC^m\to \CC^1$ has exponential type $V$
 if for any $\vep>0$ there exists a constant $C_0(\vep,f)$ such that
 for all $z\in \CC$,
 $\vert f(z)\vert\le C_0(\vep,f)\exp\left((1+\vep)\|z\|_V^*\right)$.
 \end{definition}
  The class of all entire function of exponential type $V$ is denoted
  by $B_V$.
  It is easy to verify that if $V_1\subseteq V_2$, then
  $B_{V_1}\subseteq B_{V_2}$.
  Throughout the paper, if no confusion may occur, the same notation is applied to
  $f\in B_V$ and its restriction to $\R^m$ (e.g., in the form
  $f\in  B_V\cap L_p(\R^m))$.
  The class $B_V$ was defined by Stein and Weiss
  \cite[Sect. 3.4]{SW1971}. For $V=\Pi^m(\sa),\,V=Q^m(M),$ and
  $V=\BB^m(M)$, similar
  classes were
  defined by Bernstein \cite{B1948} and Nikolskii
  \cite{N1951}, \cite[Sects. 3.1, 3.2.6]{N1969}, see also
  \cite[Definition 5.1]{DP2010}.
  Properties of functions from $B_V$ have been investigated in numerous
  publications (see, e.g., \cite{B1948, N1951, N1969, SW1971, NW1978,
  G1982, G1991, G2001} and
  references therein). Some of these properties are presented in Section \ref{S2}.
  We need more definitions for entire functions of several variables.

  We say that an entire function $f:\CC^m\to \CC^1$ has exponential type
 if
 \ba
 \limsup_{\vert z\vert\to\iy}
 \log\vert f(z)\vert\left/\sum_{j=1}^m\vert z_j\vert<\iy\right..
 \ea
 The class of all entire functions of exponential type is denoted
  by $B$.

  Next, let us define the indicator and the $P$-indicator of $f\in B$ by
  \ba
  h_f(y,x):=\limsup_{\rho\to\iy}\log\vert f(x+i\rho y)\vert,\quad
  h_f(y):=\sup_{x\in\R^m} h_f(y,x),\quad x\in\R^m,\quad y\in\R^m,
  \ea
  respectively. Note that equivalent definitions of $h_f(y,x)$ and $h_f(y)$
  were introduced by Plancherel and P\'{o}lya \cite{PP1937} (see also
  \cite[Sect. 3.4.2]{R1971}). The class of all functions $f\in B$ with
  $h_f(y)\le \| y\|_V^*,\,y\in\R^m,$ is denoted by $B^*_V$. A similar class was
  introduced by the author \cite{G1982}.
  It turns out that $B_V^*\cap L_p(\R^m)=B_V\cap L_p(\R^m),\,p\in(0,\iy]$
  (see Lemma \ref{L2.1} (d)),
  but the definition of $B_V^*$ is needed for the proof of important Lemma
  \ref{L2.2}.

Throughout the paper $C,\,C_0,\,C_1,\ldots$ denote positive constants independent
of essential parameters.
 Occasionally we indicate dependence on certain parameters.
 The same symbol $C$ does not
 necessarily denote the same constant in different occurrences.
In addition, we use the ceiling function
 $\lceil a \rceil$.
  \vspace{.12in}\\
\textbf{ Markov-Bernstein-Nikolskii Type Inequalities.}
Let
$D_N:=\sum_{\vert\al\vert=N}b_\al D^\al$
be a linear differential
operator with constant coefficients $b_\al\in\CC^1,\,\vert\al\vert=N,\,
N\in \Z^1_+$.
We assume that $D_0$
is the identity operator.

 Next, we define  sharp constants in multivariate Markov-Bernstein-Nikolskii type
inequalities for algebraic and trigonometric polynomials
and entire functions of exponential type. Let
\bna
&&M_{p,D_N,n,m,V}:=n^{-N-m/p}
\sup_{P\in\PP_{n,m}\setminus\{0\}}\frac{\vert D_N(P)(0)\vert}
{\|P\|_{L_p(V^*)}},\label{E1.2}\\
&&P_{p,D_N,a,m,V}:=a^{-N-m/p}
\sup_{T\in\TT_{aV}\setminus\{0\}}\frac{\|D_N(T)\|_{L_\iy(Q^m(\pi))}}
{\|T\|_{L_p(Q^m(\pi))}},\nonumber\\
&& E_{p,D_N,m,V}:=
\sup_{f\in (B_{V}\cap L_p(\R^m))\setminus\{0\}}\frac{\|D_N(f)\|_{L_\iy(\R^m)}}
{\|f\|_{L_p(\R^m)}}.\label{E1.3}
\ena
Here, $a>0,\,N\in\Z^1_+,\,V\subset\R^m$, and $p\in(0,\iy]$.

The purpose of this paper is to prove limit relations between $E_{p,D_N,m,V}$ and
$M_{p,D_N,n,m,V}$.
The limit relation for multivariate trigonometric polynomials
\beq\label{E1.5}
E_{p,D_N,m,V}= \lim_{a\to\iy}P_{p,D_N,a,m,V},\qquad p\in(0,\iy],
\eeq
was proved by the author \cite[Theorem 1.3]{G2018}.
In the univariate case of $V=[-1,1],\,D_N=d^N/dx^N$, and $a\in\N$,
  \eqref{E1.5} was proved by the author and Tikhonov \cite{GT2017}.
In  earlier publications \cite{LL2015a, LL2015b}, Levin and Lubinsky established
versions of \eqref{E1.5} on the unit circle for $N=0$.
Certain extensions of the Levin-Lubinsky's results to the $m$-dimensional
unit sphere in $\R^{m+1}$ were recently proved by Dai, Gorbachev, and Tikhonov
\cite{DGT2018}.

In case of an even $N\in\Z^1_+,\,p\in[1,\iy]$, the unit ball $V=\BB^m(1)$, and the operator
$D_N=\Delta^{N/2}$, where $\Delta$ is the Laplace operator,
the author \cite[Corollary 4.4]{G2019} proved the
following limit relation
for multivariate algebraic polynomials:
\beq\label{E1.6}
E_{p,D_N,m,V}= \lim_{n\to\iy}M_{p,D_N,n,m,V}.
\eeq
The proof of \eqref{E1.6} was based on invariance theorems and limit
relations for sharp constants in univariate weighted spaces
(see \cite[Theorems 2.1, 2.2, 4.1]{G2019}).
Note that certain properties of the sharp constants in
univariate weighted spaces are discussed by Arestov and Deikalova
\cite{AD2015}.
For $m=1,\,V=[-1,1],\,D_N=d^N/dx^N,\,N\in\Z^1_+$,
and $p\in(0,\iy]$, \eqref{E1.6} was established in
\cite[Theorem 1.1]{G2017}.

In this paper we extend relation \eqref{E1.6} to any $N\in\Z^1_+,\,
p\in(0,\iy],\,V\subset\R^m$, and $D_N$.
\vspace{.12in}\\
 \textbf{Main Results and Remarks.} We recall that
 $V$ is a centrally symmetric (with respect to the origin)
 closed convex body in $\R^m$.

 \begin{theorem} \label{T1.2}
 If $n\in\N,\,N\in\Z^1_+,\,V\subset\R^m$, and $p\in(0,\iy]$, then
  $ \lim_{n\to\iy}M_{p,D_N,n,m,V}$ exists and
 \beq \label{E1.7}
  E_{p,D_N,m,V}=\lim_{n\to\iy}M_{p,D_N,n,m,V}.
 \eeq
 In addition, there exists a nontrivial function $f_0\in  B_V\cap L_p(\R^m)$ such that
 \beq \label{E1.8}
 \|D_N(f_0)\|_{L_\iy(\R^m)}/\|f_0\|_{L_p(\R^m)}=\lim_{n\to\iy}M_{p,D_N,n,m,V}.
 \eeq
 \end{theorem}
 \noindent
 The following corollary is a direct consequence of Theorem \ref{T1.2} and
 relation \eqref{E1.5}.
 \begin{corollary} \label{C1.3}
 If $n\in\N,\,N\in\Z^1_+,\,V\subset\R^m$, and $p\in(0,\iy]$,
 then
 \ba
 E_{p,D_N,m,V}=\lim_{n\to\iy}M_{p,D_N,n,m,V}
 =\lim_{a\to\iy}P_{p,D_N,a,m,V}.
 \ea
 \end{corollary}

 \begin{remark}\label{R1.4}
Relations \eqref{E1.7} and \eqref{E1.8} show that the function
$f_0\in B_V\cap L_p(\R^m)$
from Theorem \ref{T1.2} is an extremal function for $E_{p,D_N,m,V}$.
 \end{remark}

\begin{remark}\label{R1.5}
In definitions \eqref{E1.2} and \eqref{E1.3} of the sharp constants we
discuss only complex-valued functions $P$ and $f$. We can define similarly
the "real" sharp constants if the suprema in \eqref{E1.2} and \eqref{E1.3}
 are taken over all real-valued functions
on $\R^m$ from $\PP_{n,m}\setminus\{0\}$
and $(B_V\cap L_p(\R^m))\setminus\{0\}$, respectively.
It turns out that the "complex" and "real" sharp constants coincide.
 For $m=1$ this fact was proved in \cite[Sect. 1]{G2017} (cf.
\cite[Theorem 1.1]{GT2017}) and the case of $m>1$ can be proved similarly.
\end{remark}

\begin{remark}\label{R1.5a}
Note that the following relation for a different class of polynomials holds true:
\beq \label{E1.8a}
 E_{p,D_N,m,Q^m(1)}
 =\lim_{n\to\iy}n^{-N-m/p}
 \sup_{P}\frac{\vert D_N(P)(0)\vert}
 {\|P\|_{L_p(Q^m(1))}},
 \eeq
where the supremum in \eqref{E1.8a} is taken over all nontrivial
polynomials $P$ in $m$ variables of degree at most $n$ in each variable.
In addition, note that, unlike relation  \eqref{E1.7}, both sides of
\eqref{E1.8a} contain the same set $V= Q^m(1)$.
\end{remark}

\begin{remark}\label{R1.6}
Equality \eqref{E1.7} presents one more example of a
relation between
entire functions of exponential type $V$ on $\R^m$ and polynomials on
$V^*$. The limit relations for the corresponding errors of approximation
of $f\in L_p(\R^m)$ by functions from $B_V$ and
polynomials from $\PP_{n,m}$ restricted to $n V^*$
 were proved by the author \cite[Theorems 5.1, 5.2]{G1982}.
  Certain versions of
 these results are discussed below in Lemmas \ref{L2.3}
 and \ref{L2.4}.
\end{remark}

\begin{remark}\label{R1.7}
Let us  introduce a sharp constant in the classic inequality
of different metrics
\ba
N_{p,D_0,n,m,\Omega}:=n^{-\mu}
\sup_{P\in\PP_{n,m}\setminus\{0\}}\frac{\| P\|_{L_\iy(\Omega)}}
{\|P\|_{L_p(\Omega)}},\qquad p\in(0,\iy),
\ea
where $\Omega$ is a compact subset of $\R^m$ and
$\mu=\mu(\Omega,p,m)>0$  is a constant. Let us assume that
the following
inequalities hold true:
\beq \label{E1.9}
0<\liminf_{n\to\iy} N_{p,D_0,n,m,\Omega}
\le \limsup_{n\to\iy} N_{p,D_0,n,m,\Omega}<\iy,
\eeq
if any.
For example, the third inequality in \eqref{E1.9}
 holds true for $\mu=2m/p$ and any domain $\Omega$,
 satisfying the cone condition, in particular for
 convex bodies (not necessarily symmetric),
 see \cite[Theorem 1]{D1972}, \cite[Theorem 2]{D1976},
 \cite[Theorem 2]{G1978}. It is also valid for
 $\mu=(m+1)/p$ and any domain $\Omega$ with the smooth boundary,
 see \cite[Theorem 2]{D1972}, \cite[Theorem 5]{D1976},
 \cite[p. 433]{KS1997}. The typical examples
 of sets $\Omega$ when all
  inequalities in \eqref{E1.9} hold true are the unit
  cube $\Omega=Q^m(1)$ for $\mu=2m/p$  and the unit ball
   $\Omega=\BB^m(1)$
  for $\mu=(m+1)/p$. More examples and discussions are
  presented by Ditzian and Prymak \cite{DP2016}.
  The author (see \cite[Theorem 1.4]{G2017} and
 \cite[Corollary 4.6]{G2019}) proved that
 $\lim_{n\to\iy} N_{p,D_0,n,1,[-1,1]}$ exists
 for $\mu=2/p$ and
 found its exact value.
  However, it is unknown as to whether
  $\lim_{n\to\iy}N_{p,D_0,n,m,\Omega}$ exists for $m>1$.
\end{remark}

The proof of Theorem \ref {T1.2} is presented in
Section \ref{S3}.
Section \ref{S2} contains certain properties of functions
from $B_V$ and $\PP_{n,m}$.
In particular, Lemma \ref{L2.7} that discusses inequalities
 of different metrics for multivariate polynomials on smaller
  domains is of independent interest.

  \section{Properties of Entire Functions and
  Polynomials}\label{S2}
\setcounter{equation}{0}
\noindent
In this section we discuss certain properties of multivariate entire functions
 of exponential type and polynomials
that are needed for the proof of Theorem \ref {T1.2}.
We start with certain properties of entire functions
 of exponential type.
 \begin{lemma}\label{L2.1}
 (a) If $f\in B_V$, then there exists $M=M(V)$ such that $f\in B_{Q^m(M)}$.\\
 (b) $f\in B$ if and only if there exists $V\subset \R^m$ such that $f\in B_V$.\\
 (c) The following Bernstein and Nikolskii type inequalities hold true:
 \bna
 &&\left\|D^\al(f)\right\|_{L_{\iy}(\R^m)}
  \le C
  \left\|f\right\|_{L_{\iy}(\R^m)},\quad f\in B_V\cap L_\iy(\R^m),
  \quad \al\in Z^m_+,\label{E2.1a}\\
   &&\left\|f\right\|_{L_{\iy}(\R^m)}
  \le C
  \left\|f\right\|_{L_{p}(\R^m)},\quad f\in B_V\cap L_p(\R^m),\quad
 p\in(0,\iy),\label{E2.1}
  \ena
  where $C$ is independent of $f$.\\
  (d) For $p\in (0,\iy],\,B_V\cap L_p(\R^m)=B_V^*\cap L_p(\R^m)$.\\
  (e) The Lebesgue measure of the set
  $\{x\in\R^m: h_f(y,x)<h_f(y)\}$ is zero for every $y\in\R^m$.
  \end{lemma}
  \proof
  Statement (a) follows from the obvious inclusion $V\subseteq Q^m(M)$
  for a certain $M=M(V)>0$. Next, statement (b) follows from
  the inequalities
  $
  C_1\sum_{j=1}^m\vert z_j\vert
  \le \|z\|_{V}^*\le C_2\sum_{j=1}^m\vert z_j\vert,
  $
  where $C_1=C_1(V):=1/(2\max\{1,\la\})$
  (here,
  $\la$ is the least number such that $Q^m(1)\subseteq \la V$)
  and
  $C_2=C_2(V):=\max_{1\le j\le m}\max_{t\in V}\vert t_j\vert$.
  Inequality \eqref{E2.1a} for $V=Q^m(M),\,M>0$, is well known (see, e.g.,
  \cite[Eq. 3.2.2(8)]{N1969}). Then \eqref{E2.1a}
  follows from statement (a).
  Inequality \eqref{E2.1} was proved in \cite[Theorem 5.7]{NW1978}.

  To prove statement (d), we note that if $f\in B_V\cap L_p(\R^m),
  \,p\in(0,\iy)$, then $f\in B_V\cap L_\iy(\R^m)$ by \eqref{E2.1},
  while the relation $B_V\cap L_\iy(\R^m)=B_V^*\cap L_\iy(\R^m)$
  was proved in a more general form by the author
  \cite[Corollary 2.1]{G1982}.
  Statement (e) was proved by Ronkin
  (see \cite{R1958} and \cite[Theorem 3.4.3]{R1971}).\hfill $\Box$

  A compactness theorem for a set of entire functions
  is discussed below.
    \begin{lemma}\label{L2.2}
     (a) Let $\E_m$ be the set  of all multivariate entire functions
     $f(z)=\sum_{k=0}^\iy Q_k(z)$,
     where $Q_k(z)=\sum_{\vert \be\vert=k}c_\be z^\be$
     is the homogeneous polynomial of degree $k\in\Z_+^1$,
     satisfying the following conditions: for any $\vep>0$
     the following inequalities are valid:
     \bna
     &&\vert Q_k(z)\vert \le \frac{C_3(\vep)
     \left(C_4(1+\vep)\|z\|^*_V\right)^{k}}{k!},\qquad k\in\Z_+^1,
     \qquad z\in\CC^m,\label{E2.2}\\
     &&\vert Q_k(y)\vert \le \frac{C_5(\vep)
     \left((1+\vep)\|y\|^*_V\right)^{k}}{k!},\qquad k\in\Z_+^1,
     \qquad y\in\R^m,\label{E2.3}
     \ena
     where $C_3(\vep)$ and $C_5(\vep)$ are
     independent of $z,\,f$ and $k$, and
     $C_4$ is
     independent of $z,\,\vep,\,f$ and $k$.
     Then for any sequence $\{f_n\}_{n=1}^\iy\subseteq\E_m$ there exist
     a subsequence $\{f_{n_s}\}_{s=1}^\iy$ and a function $f_0\in B_{V}^*$
     such that for every $\al\in\Z_+^m$,
      \beq\label{E2.4}
     \lim_{s\to\iy} D^\al(f_{n_s})=D^\al(f_0)
     \eeq
     uniformly on each compact subset of $\CC^m$.\\
     (b) If the function $f_0$ from statement (a) belongs to
     $L_\iy(\R^m)$, then $f_0\in B_V\cap L_\iy(\R^m)$.
     \end{lemma}
     \proof
     (a) By a multivariate version of  Weierstrass' theorem
     (see, e.g., \cite[Sect. 1.2.5]{S1992}),
     it suffices to prove statement (a) for $\al=0$.
     Let us set $d_R:=\{z\in\CC^m:\| z\|_V^*\le R\},\,R>0$,
     and let
     \ba
     f_n(z)=\sum_{k=0}^\iy Q_{k,n}(z)
     =\sum_{k=0}^\iy
     \sum_{\vert \be\vert=k}c_{\be,n} z^\be,\qquad n\in\N.
     \ea
     By \eqref{E2.2}, the polynomials $Q_{k,n},\, n\in\N$,
     are uniformly bounded on $d_R$ for a fixed $R\ge 0$
     and each $k\in\Z^1_+$. Therefore, using the Cantor
     diagonal process, we can select a subsequence
     $\{f_{n_s}\}_{s=1}^\iy$ such that
     \beq\label{E2.5}
     \lim_{s\to\iy}Q_{k,n_s}=Q_{k,0}
     \eeq
     uniformly on $d_R$ for any $R\ge 0$ and $k\in\Z^1_+$.
     Here, $Q_{k,0}=\sum_{\vert \be\vert=k}c_{\be,0} z^\be$
     are homogeneous polynomials of degree $k,\,k\in\Z^1_+$.
     Then inequality \eqref{E2.2} shows that
      $f_0(z)=\sum_{k=0}^\iy Q_{k,0}(z)$ is an entire function
       that belongs to the class $B_{C_4V}$.
     Further, we obtain for $M\in\N$
  \ba
 && \max_{z\in d_R}\left\vert
  f_{0}(z)-f_{n_s}(z)\right\vert\\
  &&\le \max_{z\in d_R}
  \sum_{k=0}^{M-1} \left(
  \left\vert
  Q_{k,0}(z)-Q_{k,n_s}(z)\right\vert\right)
  + \max_{z\in d_R}
  \sum_{k=M}^{\iy}\left(
  \left\vert
  Q_{k,0}(z)\vert+\vert Q_{k,n_s}(z)\right\vert\right)
  =S_1+S_2,
  \ea
   where by  \eqref{E2.2} for $\vep=1$ and by \eqref{E2.5},
   \ba
   S_2\le C_3(1)\sum_{k=M}^{\iy}\frac{(2C_4R)^{k}}{k!}.
   \ea
   Then given $\de>0$ and $R>0$, we can choose $M=M(\de,R)$
   such that $S_2<\de/2$. Finally by \eqref{E2.5}, we can
   choose $s_0=s_0(\de,R)\in\N$ such that $S_1<\de/2$ for all $s\ge s_0$.
   Thus \eqref{E2.4} holds
   uniformly on $d_R$, where $f_0\in B_{C_4V}$.
   In addition, $f_0\in B$ by Lemma \ref{L2.1} (b).

    Next, we prove that $f_0\in B^*_V$. It follows from \eqref{E2.3} and
    \eqref{E2.5} that for any $\vep>0$ and $y\in\R^m$,
    \beq\label{E2.6}
    \vert f_0(iy)\vert\le \sum_{k=0}^\iy\vert Q_{k,0}(y)\vert
    \le C_5(\vep)\exp((1+\vep)\|y\|^*_V).
    \eeq
    The function $f_0(z)\exp(-(1+\vep)\|z\|^*_V)$ is continuous
     at $z=iy,\,y\in\R^m$. Hence by \eqref{E2.6}, there exists a number
     $\vep_1=\vep_1(y)$ such that for any
     $x\in R^m(\vep_1),\,\vert f_0(x+iy)
     \vert\le 2C_5(\vep)\exp((1+\vep)\|y\|^*_V)$.
     Therefore,  for a fixed $y\in\R^m$ and any $x\in R^m(\vep_1),$
     \beq\label{E2.7}
      h_{f_0}(y,x)=\limsup_{\rho\to\iy}\log\vert f_0(x+i\rho y)\vert
      \le (1+\vep)\|y\|^*_V.
     \eeq
     By Lemma \ref{L2.1} (e), the Lebesgue measure of the set
     $\{x\in\R^m: h_{f_0}(y,x)<h_{f_0}(y)\}$ is equal to zero for every fixed $y\in\R^m$.
     Hence using \eqref{E2.7}, we see that there
      exists $x_0=x_0(y)\in R^m(\vep_1)$ such that
     \ba
      h_{f_0}(y)=\sup_{x\in\R^m} h_{f_0}(y,x)= h_{f_0}(y,x_0)
      \le (1+\vep)\|y\|^*_V.
      \ea
      Thus $f_0\in\bigcap_{\de>0}B_{(1+\vep)V}^*=B_{V}^*$.
      This completes the proof of statement (a) of the lemma.\\
    (b) If $f_0\in L_\iy(\R^m)$, then $f_0\in B_V^*\cap L_\iy(\R^m)$
      by statement (a), and $f_0\in B_V\cap L_\iy(\R^m)$
      by Lemma \ref{L2.1} (d).
   \hfill $\Box$

In the next two lemmas we discuss estimates of the error of polynomial approximation
  for functions from $B_V$.

\begin{lemma}\label{L2.3}
  For any function $f\in B_V\cap L_\iy(\R^m),\,\tau\in(0,1)$, and $k\in\N$,
   there is a polynomial
  $P_k\in\PP_{k,m}$ such that
  \beq\label{E2.8}
  \|f-P_k\|_{L_\iy(\tau kV^*)}
  \le C_6(\tau,m,V)k^{\frac{2m}{m+4}}\exp(-C_7(\tau,m)\, k)\|f\|_{L_\iy(\R^m)}.
  \eeq
  \end{lemma}
  \noindent
  This result was proved by the author \cite[Lemma 4.4]{G1982}.
  In case of $m=1$ a version of Lemma \ref{L2.3} was established
  by Bernstein \cite{B1946}
  (see also \cite[Sect. 5.4.4]{T1963} and
  \cite[Appendix, Sect. 83]{A1965}).  More general and
  more precise inequalities were obtained
  in \cite{G1982} and \cite{G1991}.

  \begin{lemma}\label{L2.4}
  For any $f\in  B_V\cap L_\iy(\R^m),\,\tau\in(0,1)$, and $n\in\N$,
  there is a polynomial
  $P_n\in\PP_{n,m}$ such that for $\al\in\Z_+^m$ and
    $r\in(0,\iy]$,
   \beq\label{E2.9}
  \lim_{n\to\iy}\left\|D^\al(f)-D^\al(P_n)\right\|_{L_{r}(\tau nV^*)}=0.
  \eeq
  \end{lemma}
  \proof
  First of all, for $P_k\in\PP_{k,m},\,k\in Z_+^1,\,a>0$, and $\al\in\Z_+^m$,
   we need the following crude Markov-type inequality:
  \beq\label{E2.10}
\left\|D^\al(P_k)\right\|_{L_\iy(aV^*)}
  \le \left(\frac{4k^2}{a\, w(V^*)}
  \right)^{\vert\al\vert}\|P_k\|_{L_\iy(aV^*)}.
  \eeq
  Inequality \eqref{E2.10} easily follows from a multivariate
  A. A. Markov-type
  inequality proved by Wilhelmsen \cite[Theorem 3.1]{W1974}.

  Next, let $\{P_k\}_{k=1}^\iy$ be the sequence of
  polynomials from Lemma \ref{L2.3}. Then using \eqref{E2.10}
  for $a=\tau k$
  and  \eqref{E2.8}, we obtain
  \ba
  &&\left\|D^\al(f)-D^\al(P_n)
  \right\|_{L_\iy(\tau n V^*)}
  \le \sum_{k=n}^\iy \left\|D^\al(P_k-P_{k+1})
  \right\|_{L_\iy(\tau n V^*)}\\
  &&\le (4/w(V^*))^{\vert\al\vert}
  (\tau n)^{-\vert\al\vert} \sum_{k=n}^\iy
  (k+1)^{2\vert\al\vert}
      \left\|P_k-P_{k+1}\right\|_{L_\iy(\tau n V^*)}\\
  &&\le  (4/w(V^*))^{\vert\al\vert}(\tau n)^{-\vert\al\vert}
  \sum_{k=n}^\iy (k+1)^{2\vert\al\vert}
      \left(\left\|f-P_{k}\right\|_{L_\iy(\tau k V^*)}
 + \left\|f-P_{k+1}\right\|_{L_\iy(\tau (k+1) V^*)}\right)\\
  &&\le  2C_6(4/w(V^*))^{\vert\al\vert}(\tau n)^{-\vert\al\vert}
  \sum_{k=n}^\iy (k+1)^{2\vert\al\vert+2}\exp(-C_7 k)\,\|f\|_{L_\iy(\R^m)}.
  \ea
  Hence for $n\in\N,\,\al\in\Z_+^m$, and
    $r\in(0,\iy]$
    we have
  \ba
  &&\left\|D^\al(f)-D^\al(P_n)
  \right\|_{L_{r}(\tau n V^*)}\\
  &&\le \vert V^*\vert_m^{1/r} (\tau n)^{m/r}
   \left\|D^\al(f)
  -D^\al(P_n)
  \right\|_{L_\iy(\tau n V^*)}\\
  &&\le C_8(\tau,r,m,V,\vert\al\vert) n^{m/r-\vert\al\vert}
      \int_{n}^\iy \la^{2\vert\al\vert+2}\exp(-C_7 \la)d\la\, \|f\|_{L_\iy(\R^m)}\\
   &&\le C_9(\tau,r,m,V,\vert\al\vert)
   n^{m/r+\vert\al\vert+2}\exp(-C_7 n)\,\|f\|_{L_\iy(\R^m)}.
  \ea
  Thus \eqref{E2.9} is established.
  \hfill $\Box$

 Certain inequalities for multivariate
 trigonometric and algebraic polynomials
     are discussed in the next three lemmas.
  \begin{lemma}\label{L2.5}
  For a polynomial
  $P(y)=\sum_{\vert\be\vert\le n}c_\be y^\be \in\PP_{n,m},\,n\in\Z^1_+$,
  the following inequalities hold true:
  \bna
  &&\left\vert \sum_{\vert\be\vert = k}c_\be z^\be\right\vert
  \le \frac{\left( C_4 n\|z\|^*_V\right)^k}{k! a^k}
  \|P\|_{L_\iy(aV)},\quad z\in\CC^m,\quad 0\le k\le n,\quad a>0,
  \label{E2.11a}\\
  &&\left\vert \sum_{\vert\be\vert = k}c_\be y^\be\right\vert
  \le \frac{\left( n\|y\|^*_V\right)^k}{k! a^k}
  \|P\|_{L_\iy(aV)},\quad y\in\R^m,\quad 0\le k\le n,\quad a>0,
  \label{E2.11}
  \ena
  where $C_4$ is independent of $z,\,a,\,f$, and $k$.
  \end{lemma}
  \proof
  Inequality \eqref{E2.11a} is proved in \cite[Lemma 3.2]{G1982}.
  We first prove inequality \eqref{E2.11} for $m=1$ and $V=[-1,1]$,
  i.e., we prove the inequality
  \beq\label{E2.12}
  \left\vert d_k\tau^k\right\vert
  \le \frac{\left(n\vert\tau\vert\right)^k}{k! a^k}
  \max_{\tau\in[-a,a]}\left\vert \sum_{k=0}^n d_k\tau^k\right\vert,
  \quad \tau\in\R^1,\quad 0\le k\le n,
  \eeq
  which is equivalent to
   \beq\label{E2.13}
   \vert d_k\vert
   \le \frac{M \,n^k }{k! a^k},\qquad 0\le k\le n,
   \eeq
   where $M:=
   \max_{\tau\in[-a,a]}\left\vert \sum_{k=0}^n d_k\tau^k\right\vert$.
   This inequality without proof was discussed in
   \cite{B1946} and \cite[Eq. 2.6(9)]{T1963} as a
   corollary of V. A. Markov's inequality \cite{M1892} (see also
   \cite[Sect. 2.9.12]{T1963}) and
   \cite[Eqs. (5.1.4.1) and (6.1.2.5)]{MMR1994})
   \beq\label{E2.14}
    \vert d_k\vert \le
    \left\{\begin{array}{ll}
 \frac{M\left\vert T_{n-1}^{(k)}(0)\right\vert}{k! a^k}, &n-k\,\mbox{is odd},\\
  \frac{M\left\vert T_{n}^{(k)}(0)\right\vert}{k! a^k}, &n-k\,\mbox{is even},
 \end{array}\right.
   \eeq
   where $T_p\in\PP_{p,1}$ is the Chebyshev polynomial of the first kind.
   It is sufficient to derive \eqref{E2.13} from \eqref{E2.14}
   for $M=1$ and $a=1$. We consider three cases.\\
   \textbf{Case 1:} $k=2p,\,n=2N$ or $n=2N+1,\,0\le p\le N$. Then we have from
   \eqref{E2.14}
   \bna\label{E2.15}
    \vert d_{2p}\vert \le \frac{2^{2p}N}{N+p}\binom{N+p}{2p}
   = \left\{\begin{array}{ll}
    \frac{2^{2p}\prod_{l=0}^{p-1}(N^2-l^2)}{(2p)!}, &1\le p\le N,\\
    1,&p=0 \end{array}\right.
    \le \frac{(2N)^{2p}}{(2p)!}\le \frac{n^{2p}}{(2p)!}.
    \ena
    \textbf{Case 2:} $k=2p+1,\,n=2N+1,\,0\le p\le N$. Then we have from
   \eqref{E2.14}
   \bna\label{E2.16}
    \vert d_{2p+1}\vert &\le& \frac{2^{2p}(2N+1)}{2p+1}\binom{N+p}{2p}
   = \left\{\begin{array}{ll}
    \frac{2^{2p}(2N+1)\prod_{l=1}^{p}(N^2-l^2+N+l)}{(2p+1)!}, &1\le p\le N,\\
    2N+1,&p=0 \end{array}\right.\nonumber\\
    &\le& \frac{(2N+1)^{2p+1}}{(2p+1)!}= \frac{n^{2p+1}}{(2p+1)!}.
    \ena
    \textbf{Case 3:} $k=2p+1,\,n=2N,\,0\le p\le N-1$. It follows from
    \eqref{E2.14} that we can replace $N$ with $N-1$ in inequalities
    \eqref{E2.16}. Therefore,
    \beq\label{E2.17}
   \vert d_{2p+1}\vert
   \le \frac{(2N-1)^{2p+1}}{(2p+1)!}\le \frac{n^{2p+1}}{(2p+1)!}.
   \eeq
   Thus \eqref{E2.13} and \eqref{E2.12} follow from \eqref{E2.15},
   \eqref{E2.16}, and \eqref{E2.17}.

   To prove \eqref{E2.11} for $m>1$, we first draw a line $l_y$ passing through
   a fixed point $y\in\R^m$ and the origin with the equation
   $t=\tau v_y$, where $\tau\in\R^1$ and
   $v_y\in l_y\cap\partial(aV^*)$
   is the closest point to $y$.
   In particular, $\vert\tau\vert\le 1$ for $t\in l_y\cap aV^*$.
   In addition,  we see by \eqref{E1.1a} that for $t=y$,
   \beq\label{E2.18}
   \vert\tau\vert= \vert\tau(y)\vert=\vert y\vert/\vert v_y\vert
   =\|y\|^*_V/a.
   \eeq
   Next, the restriction of $P\in\PP_{n,m}$ to $l_y$ is a polynomial
   $p(\tau)=\sum_{k=0}^nd_k\tau^k\in\PP_{n,1}$, where
   $d_k=\sum_{\vert\be\vert = k}c_\be v_y^\be,\,0\le k\le n$,
    and using relations
    \eqref{E2.12} and \eqref{E2.18}, we obtain
    \ba
    \left\vert \sum_{\vert\be\vert = k}c_\be y^\be\right\vert
    =\vert d_k\vert \vert\tau\vert^k
    \le \frac{(n\vert\tau\vert)^k}
    {k!}\|p\|_{L_\iy(l_y\cap aV^*)}
    \le \frac{(n\|y\|^*_V)^k}
    {k!a^k}\|P\|_{L_\iy(aV^*)},
    \quad 0\le k\le n.
    \ea
    Thus \eqref{E2.11} is established.
  \hfill $\Box$

  Note that a version of Lemma \ref{L2.5} was proved in
  \cite[Theorem 1]{G2002}.

  Next, we need a multivariate version of a generalized Bari's
  inequality \cite[Lemma 2.4]{G2017}. Bari \cite[Theorem 6]{Ba1954}
  proved the following result  for $m=1$ and $p\in[1,\iy)$.

  \begin{lemma}\label{L2.6}
  Let $T\in\TT_{nQ^m(1)}$ be an even trigonometric polynomial
  of degree at most $n$ in each variable, $n\in\N$,
   and let $0\le a_j<a_{1,j}
  \le b_{1,j}<b_j\le\pi,\,1\le j\le m$.
  Then for  $p\in(0,\iy)$,
  \beq\label{E2.19}
  \|T\|_{L_\iy\left(\prod_{j=1}^m[a_{1,j},b_{1,j}]\right)}
  \le C_{10}n^{m/p}\|T\|_{L_p\left(\prod_{j=1}^m[a_j,b_j]\right)}.
  \eeq
  where $C_{10}$ is independent of $n$ and $T$.
  \end{lemma}
  \proof
  The proof follows that of \cite[Theorem 6]{Ba1954}
  and \cite[Lemma 2.4]{G2017}. Let
  $t_j=t_j(u_j):[a_j,b_j]\to[0,\pi]$ be the one-to-one function defined by the equation
  \beq\label{E2.20}
  \cos u_j=\frac{\cos a_j-\cos b_j}{2}\cos t_j+\frac{\cos a_j+\cos b_j}{2},\qquad 1\le j\le m,
  \eeq
  and let $\la_{1,j}=t_j(a_{1,j}),\,\la_{2,j}=t_j(b_{1,j})$,
  where $0<\la_{1,j}\le \la_{2,j}<\pi,\,1\le j\le m$.
  Then
  \beq\label{E2.21}
  \sin u_j\,du_j=\frac{\cos a_j-\cos b_j}{2}\sin t_j\,dt_j,\qquad 1\le j\le m,
  \eeq
   and
   \beq\label{E2.22}
   \de:=\min_{1\le j\le m}\inf_{t_j\in[\la_{1,j},\la_{2,j}]}\sin t_j
   =\min_{1\le j\le m}\min\{\sin \la_{1,j},\sin \la_{2,j}\}>0.
   \eeq
   Next, let $H_n(t):=T(u)$, by substitution \eqref{E2.20}, and let
   $G_{n+d}(t):=H_n(t)\prod_{j=1}^m\sin^d t_j\in\TT_{(n+d)Q^m(1)}$,
    where $d:=\lceil 1/p\rceil$, i.e., $dp\ge 1$.  Then using
    \eqref{E2.21}, \eqref{E2.22}, and Nikolskii's inequality for
    multivariate trigonometric polynomials of degree at most $n+d$
    in each variable \cite[Sect. 3.3.3]{N1969},
    we obtain
    \ba
     &&\|T\|_{L_\iy\left(\prod_{j=1}^m[a_{1,j},b_{1,j}]\right)}
     =\|H_n\|_{L_\iy\left(\prod_{j=1}^m[\la_{1,j},\la_{2,j}]\right)}
     \le\de^{-md}\|G_{n+d}\|
     _{L_\iy\left(\prod_{j=1}^m[\la_{1,j},\la_{2,j}]\right)}\\
     &&\le \de^{-md}\|G_{n+d}\|
     _{L_\iy\left([-\pi,\pi]^m\right)}
     \le \de^{-md}C_{11}(p,m)(n+d)^{m/p}\|G_{n+d}\|
     _{L_p\left([-\pi,\pi]^m\right)}\\
     &&= 2^{m/p}\de^{-md}C_{11}(p,m)(n+d)^{m/p}
     \left(\int_{[0,\pi]^m}
     \vert H_n(t)\vert^p\prod_{j=1}^m\sin^{dp}t_j\,dt\right)^{1/p}\\
     &&\le 2^{m/p}\de^{-md}C_{11}(p,m)(n+d)^{m/p}
     \left(\int_{[0,\pi]^m}
     \vert H_n(t)\vert^p\prod_{j=1}^m\sin t_j\,dt\right)^{1/p}\\
     &&= \frac{2^{2m/p}C_{11}(p,m)(n+d)^{m/p}}
     {\de^{md}\prod_{j=1}^m(\cos a_j-\cos b_j)^{m/p}}
     \left(\int_{\prod_{j=1}^m[a_j,b_j]}
      \vert T(u)\vert^p\prod_{j=1}^m\sin u_j\,du\right)^{1/p}\\
     &&\le C_{10} n^{m/p}\|T\|_{L_p(\prod_{j=1}^m[a_j,b_j])}.
     \ea
     Hence \eqref{E2.19} follows.\hfill $\Box$\vspace{.1in}\\
      Certain polynomial estimates based on Lemma \ref{L2.6}
      are discussed in the following lemma.
    \begin{lemma}\label{L2.7}
    Let $P$ be an algebraic polynomial in $m$ variables of degree
    at most $n$ in each variable.\\
    (a) For  $p\in(0,\iy)$ and $M>0$,
    \beq\label{E2.23}
    \vert P(0)\vert\le C_{12}(p,m)(n/M)^{m/p}\|P\|_{L_p(Q^m(M))}.
    \eeq
    (b) For  $p\in(0,\iy),\,a>0$, and $\vep>0$, the
    following inequality holds true:
    \beq\label{E2.24}
    \|P\|_{L_\iy(aV^*)}
    \le C_{13}(\vep,p,m,V) (n/a)^{m/p}
    \|P\|_{L_p((1+\vep)aV^*)}.
    \eeq
    \end{lemma}
    \proof
    Setting
    \ba
    T(t):=P(M\cos t_1,\ldots, M\cos t_m),\,\,
    a_j=\pi/3,\,\, b_j=2\pi/3,\,\, a_{1,j}=b_{1.j}=\pi/2,\,\, 1\le j\le m,
   \ea
     we have by Lemma \ref{L2.6}
     \ba
     \vert P(0)\vert= \vert T(\pi/2,\ldots,\pi/2)\vert
     &\le & C_{14}(p,m)n^{m/p}
     \left(\int_{Q^m(1/2)}\vert P(Mx)\vert^p
     \prod_{j=1}^m(1-x_j^2)^{-1/2}dx\right)^{1/p}\nonumber\\
     &\le &C_{12}(p,m) (n/M)^{m/p}\|P\|_{L_p(Q^m(M))}.
     \ea
     Hence \eqref{E2.23} holds true. \\
     (b) Let $y_0\in aV^*$ such that
     $\|P\|_{L_\iy(aV^*)}=\left\vert P(y_0)\right\vert$.
     The cube $Q_M:=Q^m(aM)+y_0$, where $M:=\vep\, a\, w(V^*)/(2\sqrt{m})$,
     is a subset of $(1+\vep)aV^*$. Indeed, for any $y\in Q_M$
     we have by \eqref{E1.1a},
     \ba
     \|y\|_V^*\le \|y_0\|_V^*+\|y-y_0\|_V^*
     \le a +\frac{\vert y-y_0\vert}{\vert v_{y-y_0}\vert}
     \le a +\frac{2\vert y-y_0\vert}{w(V^*)}
     \le a(1+\vep).
     \ea
      Therefore, it follows from
     \eqref{E2.23} that
     \ba
     \left\vert P(y_0)\right\vert
     \le C_{12}(p,m) (n/M)^{m/p}\|P\|_{L_p(Q_M)}
     \le C_{13}(\vep,p,m,V) (n/a)^{m/p}\|P\|_{L_p((1+\vep)aV^*)}.
     \ea
     This proves \eqref{E2.24}.
     \hfill $\Box$ \begin{remark}\label{R2.8}
     Lemma \ref{L2.7} (a) for $p\in[1,\iy)$ was proved
     by a different method in \cite[Lemma 3.3]{G1982}.
     Note that  \eqref{E2.24} is valid for any convex body $V^*$.
     In addition, note that an estimate of the constant in
     Nikolskii-type inequality
     \eqref{E2.24} on a smaller convex body is $C(n/a)^{m/p}$,
     while for the corresponding inequality on the same domain
     the estimate is  $C(n/a)^{(m+1)/p}$ for domains with smooth
     boundaries and the estimate is
       $C(n/a)^{2m/p}$ for domains with boundaries,
     satisfying the cone condition (see
     Remark \ref{R1.7}).
     \end{remark}

\section{Proof of Theorem  \ref{T1.2}}\label{S3}
 \noindent
\setcounter{equation}{0}
Throughout the section we use the notation $\tilde{p}=\min\{1,p\},\,
p\in(0,\iy]$, introduced in Section \ref{S1}.\vspace{.1in}\\
\emph{Proof of Theorem \ref{T1.2}.}
We first prove the inequality
\beq \label{E3.1}
  E_{p,\iy,D_N,m,V}\le\liminf_{n\to\iy}M_{p,D_N,n,m,V},\qquad p\in(0,\iy].
 \eeq
Let $f$ be any function from $B_V\cap L_p(\R^m),\,p\in(0,\iy]$.
Then $f\in B_{Q^m(M)}$ by Lemma \ref{L2.1} (a); hence
$D_N(f)\in B_{Q^m(M)}$ by \cite[Sect. 3.1]{N1969}
(see also \cite[Lemma 2.1 (d)]{G2018}).
In addition,
$D_N(f)\in L_p(\R^m)$ by Bernstein's and
Nikolskii's inequalities \eqref{E2.1a} and \eqref{E2.1} and by the "triangle" inequality
\eqref{E1.1}. Therefore,
\beq\label{E3.2}
\lim_{\vert x\vert\to\iy}D_N(f)(x)=0,\qquad p\in(0,\iy).
\eeq
Indeed, since $D_N(f)\in B_{Q^m(M)}\cap L_p(\R^m)$, \eqref{E3.2} is known for $p\in[1,\iy)$
(see, e.g., \cite[Theorem 3.2.5]{N1969}), and for $p\in(0,1)$ it follows from
Nikolskii's inequality \eqref{E2.1}, since if $D_N(f)\in L_p(\R^m),\,p\in(0,1)$,
then $D_N(f)\in L_1(\R^m)$.

Let us first prove \eqref{E3.1} for $p\in(0,\iy)$. Then by \eqref{E3.2}, there exists $x_0\in\R^m$
such that $\|D_N(f)\|_{L_\iy(\R^m)}=\left\vert D_N(f)(x_0)\right\vert$.
Without loss of generality we can assume that $x_0=0$. Let $\tau\in(0,1)$ be a fixed number.
Then using polynomials $P_n\in\PP_{n,m},\,n\in\N$, from Lemma \ref{L2.4}, we obtain
for $r=\iy$ by \eqref{E1.2},
\bna\label{E3.3}
 &&\|D_N(f)\|_{L_\iy(\R^m)}=\left\vert D_N(f)(0)\right\vert\nonumber\\
 &&\le  \lim_{n\to\iy}\left\vert D_N(f)(0)-D_N(P_n)(0)\right\vert
 +\liminf_{n\to\iy}\left\vert D_N(P_n)(0)\right\vert\nonumber\\
 &&=\liminf_{n\to\iy}\left\vert D_N(P_n)(0)\right\vert
 \le \tau^{-(N+m/p)}\liminf_{n\to\iy}\left(M_{p,D_N,n,m,V}
 \left\| P_n\right\|_{L_p(\tau n V^*)}\right).
 \ena
Next, note that $f\in L_\iy(\R^m)$, by Nikolskii's inequality \eqref{E2.2}.
 Using again Lemma \ref{L2.4} (for $\al=0$ and $r=p$), we have from \eqref{E1.1}
 \bna\label{E3.4}
 \limsup_{n\to\iy}
 \left\| P_n\right\|_{L_p(\tau n V^*)}
 \le \lim_{n\to\iy}\left(\|f-P_n\|_{L_p(\tau n V^*)}^{\tilde{p}}
 +\|f\|_{L_p(\tau n V^*)}^{\tilde{p}}\right)^{1/\tilde{p}}
 =\|f\|_{L_p(\R^m)}.
 \ena
 Combining \eqref{E3.3} with \eqref{E3.4}, and letting $\tau\to 1-$,
 we arrive at \eqref{E3.1} for $p\in(0,\iy)$.

 In the case $p=\iy$, for any $\vep>0$
 there exists $x_0\in\R^m$
such that $\|D_N(f)\|_{L_\iy(\R^m)}<(1+\vep)\left\vert D_N(f)(x_0)\right\vert$.
Without loss of generality we can assume that $x_0=0$. Then similarly to
\eqref{E3.3} and \eqref{E3.4} we can obtain the inequality
\beq\label{E3.5}
\|D_N(f)\|_{L_\iy(\R^m)}
<(1+\vep)\tau^{-N}\liminf_{n\to\iy}M_{\iy,D_N,n,m,V}
 \|f\|_{L_\iy(\R^m)}.
 \eeq
 Finally letting $\tau\to 1-$ and $\vep\to 0+$ in \eqref{E3.5},
 we arrive at \eqref{E3.1} for $p=\iy$.
  This completes the proof of \eqref{E3.1}.

  Further, we will prove the
 inequality
 \beq\label{E3.6}
\limsup_{n\to\iy}M_{p,D_N,n,m,V}\le E_{p,\iy,D_N,m,V},\qquad p\in(0,\iy],
\eeq
  by constructing a nontrivial function $f_0\in B_V\cap L_p(\R^m)$,
   such that
 \beq \label{E3.7}
  \limsup_{n\to\iy}M_{p,D_N,n,m,V}
  \le\|D_N(f_0)\|_{L_\iy(\R^m)}/
\|f_0\|_{L_p(\R^m)}
\le E_{p,\iy,D_N,m,V}.
 \eeq
 Then inequalities \eqref{E3.1} and \eqref{E3.6} imply  \eqref{E1.7}.
 In addition, $f_0$ is an extremal function in \eqref{E1.7},
that is, \eqref{E1.8} is valid.

It remains to construct a nontrivial function $f_0$,
satisfying \eqref{E3.7}.
We first note that
\beq \label{E3.8}
\inf_{n\in\N}M_{p,D_N,n,m,V}\ge C_{14}(p,N,D_N,m,V).
\eeq
This inequality follows immediately from  \eqref{E3.1}.
Let $P_n\in\PP_{n,m}$ be a polynomial, satisfying the equality
\beq \label{E3.9}
M_{p,D_N,n,m,V}=n^{-N-m/p}\left\vert D_N(P_n)(0)\right\vert
/\|P_n\|_{L_p(V^*)}.
\eeq
The existence of an extremal polynomial $P_n$ in \eqref{E3.9}
can be proved by the standard compactness argument (cf. \cite{GT2017}).
Next, setting $Q_n(y):=P_n(y/n)=\sum_{\vert\be\vert\le n} c_\be y^\be$,
 we have from \eqref{E3.9} that
\beq \label{E3.10}
M_{p,D_N,n,m,V}=\left\vert D_N(Q_n)(0)\right\vert
/\|Q_n\|_{L_p(nV^*)}.
\eeq
We can assume that
\beq \label{E3.11}
\left\vert D_N(Q_n)(0)\right\vert=1.
\eeq
Then it follows from \eqref{E3.10}, \eqref{E3.11}, and \eqref{E3.8}
that
\beq \label{E3.12}
\|Q_n\|_{L_p(nV^*)}=1/M_{p,D_N,n,m,V}\le 1/C_{14}(p,N,D_N,m,V).
\eeq
Further,
for any $\vep > 0$ and $0\le k\le n$, the following
inequalities hold true:
\bna
\left\vert \sum_{\vert\be\vert=k} c_\be z^\be\right\vert
&\le& \frac{C_{13}(1+\vep)^{k+m/p}
\left(C_4\|z\|_V^*\right)^k}{k!}
 \|Q_n\|_{L_p(nV^*)}\nonumber\\
&\le& \frac{\left(C_{13}/C_{14}\right)
(1+\vep)^{k+m/p}\left(C_4\|z\|_V^*\right)^k}{k!},
\qquad z\in\CC^m, \label{E3.13}\\
\left\vert \sum_{\vert\be\vert=k} c_\be y^\be\right\vert
&\le& \frac{C_{13}(1+\vep)^{k+m/p}\left(\|y\|_V^*\right)^k}{k!}
 \|Q_n\|_{L_p(nV^*)}\nonumber\\
&\le& \frac{\left(C_{13}/C_{14}\right)
(1+\vep)^{k+m/p}\left(\|y\|_V^*\right)^k}{k!},
\qquad y\in\R^m.\label{E3.13a}
\ena
Inequalities \eqref{E3.13} and \eqref{E3.13a}  follow from
\eqref{E3.12} and Lemmas \ref{L2.5} and \ref{L2.7} (b)
for $a=n/(1+\vep)$.
Let $\{n_r\}_{r=1}^\iy$ be a subsequence of $\N$ such that
\beq \label{E3.14}
\limsup_{n\to\iy}M_{p,D_N,n,m,V}=\lim_{r\to\iy}M_{p,D_N,n_r,m,V}.
\eeq
Inequalities \eqref{E3.13} and \eqref{E3.13a}
 show that the polynomial sequence
$\{Q_{n_r}\}_{r=1}^\iy\subseteq \E_m$ satisfies the conditions
of Lemma \ref{L2.2}. Therefore, there exist
a subsequence
 $\{Q_{n_{r_s}}\}_{s=1}^\iy$ and
 a function $f_0\in B_V$ such that
 \beq \label{E3.15}
  \lim_{s\to\iy}Q_{n_{r_s}}(x)=f_0(x),\qquad
 \lim_{s\to\iy}D_N\left(Q_{n_{r_s}}\right)(x)=D_N(f_0)(x),
 \eeq
 uniformly on any cube $Q^m(M),\,M>0$. Moreover, by
  \eqref{E3.11} and \eqref{E3.15},
  \beq \label{E3.16}
  \left\vert D_N(f_{0})(0)\right\vert=1.
\eeq
In addition,
using \eqref{E1.1}, \eqref{E3.15}, \eqref{E3.10},
\eqref{E3.11}, and \eqref{E3.14},
we obtain for any cube $Q^m(M),\,M>0$,
\bna \label{E3.18}
&&\|f_0\|_{L_p(Q^m(M))}
\le \lim_{s\to\iy}
\left(\|f_0-Q_{n_{r_s}}\|_{L_p(Q^m(M))}^{\tilde{p}}
+\|Q_{n_{r_s}}\|_{L_p(Q^m(M))}^{\tilde{p}}
\right)^{1/\tilde{p}}\nonumber\\
&&\le\lim_{s\to\iy}\|Q_{n_{r_s}}\|_{L_p\left(n_{r_s}V^*\right)}
=1/
\lim_{s\to\iy}M_{p,D_N,n_{r_s},m,V}.
\ena
Next using \eqref{E3.18} and \eqref{E3.8}, we see that
\beq \label{E3.19}
\|f_0\|_{L_p(\R^m)}\le 1/C_{14}(p,N,D_N,m,V).
\eeq
Therefore, $f_0$ is a nontrivial function from $B_V\cap L_p(\R^m)$,
by \eqref{E3.16} and \eqref{E3.19}. Thus for any cube $Q^m(M),\,M>0$,
we obtain from \eqref{E3.14}, \eqref{E3.15},
 and \eqref{E3.16}
\bna \label{E3.20}
\limsup_{n\to\iy}M_{p,D_N,n,m,V}
&=&\lim_{s\to\iy}\left(\|Q_{n_{r_s}}\|
_{L_p\left(n_{r_s}V^*\right)}\right)^{-1}\nonumber\\
&\le& \lim_{s\to\iy}\left(\|Q_{n_{r_s}}\|
_{L_p(Q^m(M))}\right)^{-1}\nonumber\\
&=&\left\vert D_N(f_0)(0)\right\vert/
\|f_0\|_{L_p(Q^m(M))}\nonumber\\
&\le& \left\| D_N(f_0)\right\|_{L_\iy(\R^m)}/
\|f_0\|_{L_p(Q^m(M))}.
\ena
Finally, letting $M\to \iy$ in \eqref{E3.20}, we arrive at \eqref{E3.7}.
\hfill$\Box$ \vspace{.1in}\\
\textbf{Acknowledgements.} We are grateful to the anonymous referee
 for valuable suggestions.\vspace{.1in}\\

\end{document}